\documentclass[12pt]{amsart}
\usepackage{amsmath,amssymb,graphicx,subfigure,psfrag,pifont}
\newcommand{\NS}{MR96c:20066}
\newcommand{\Epstein}{MR93i:20036}

\newcommand{\Cannon}{MR88a:20049}
\newcommand{\Thiel}{MR95e:20052}

\newcommand{\Ethesis}{Ethesis}

\newcommand{\DShap}{MR92d:20055}

\newcommand{\Gersten}{MR96k:20073}
\newcommand{\Wang}{Wangthesis}
\newcommand{\Rebbechi}{Rebbechi}

\newcommand{\ac}{almost convexity}
\newcommand{\AC}{almost convex}

\newcommand{\CG}{$\Gamma_X(G)$}
\newcommand{\cg}{Cayley graph}
\newcommand{\cc}{Cayley complex}

\newcommand{\EM}{Eilenberg MacLane}

\newcommand{\fftp}{falsification by fellow traveler property}

\newcommand{\fn}{F$_n$}
\newcommand{\ft}{F$_3$}
\newcommand{\ftw}{F$_2$}
\newcommand{\fo}{F$_1$}

\newcommand{\fp}{finitely presented}
\newcommand{\fg}{finitely generated}

\newcommand{\ga}{$\gamma $}
\newcommand{\gam}{\gamma}

\newcommand{\gset}{generating set}

\newcommand{\iif}{isoperimetric function}

\newtheorem{thm}{Theorem}
\newtheorem{lem}{Lemma}
\newtheorem{propn}{Proposition}
\newtheorem{defn}{Definition}
\newtheorem*{cor}{Corollary}

\newtheorem*{cpfftp}{Corollary to the Proof}

\begin{document}

\title[Falsification by fellow traveler property]
  {Finiteness and the falsification by fellow traveler property}

\author[M. Elder]{Murray J. Elder}
	
\address{(As at Nov 2006)Dept. of Mathematical Sciences\\
	Stevens Institute of Technology\\
	Hoboken NJ 07030 USA}
\email{melder@stevens.edu}

\date{Published in Geometriae Dedeicata Vol 95 (2002) 103--113}

\begin{abstract}
We prove that groups enjoying the \fftp\ are of type \ft, and have at most an exponential second order \iif.
\end{abstract}


\maketitle

%%%%%%%%%%%%%%%%%%%%
\section{Introduction}
%%%%%%%%%%%%%%%%%%%%
The intriguing properties of  \ac\ and the \fftp\ have been introduced by 
\cite{\Cannon} and \cite{\NS} respectively.
It has been shown that both properties are \gset\ dependent \cite{\Thiel}, \cite{\NS};
 that groups enjoying the \fftp\ are \AC;
that  both imply a finite presentation;  that \ac\ groups have at most an
 exponential \iif\ and groups with the \fftp\ at most a quadratic 
 \iif.
The classes of groups that enjoy these properties include hyperbolic groups, virtually abelian groups \cite{\NS}, and Coxeter groups \cite{\DShap}, \cite{\Rebbechi}. 
%%\cite{DS}.
In this paper we are interested in  higher dimensional finiteness for groups with these
properties.

A group is said to be of type \fn\ if it has an \EM\ space with finite $n$-skeleton.
It follows that a group is \fg\ if and only if it is of type \fo\ and \fp\ if and only if it is
 of type \ftw, so \fn\ is a natural generalization of these finiteness properties.

Gersten showed that asynchronously combable groups with departure function
(hence asynchronously automatic groups) are of type \ft, by constructing the universal cover of
 an \EM\ space  having finitely many types of 3-cells \cite{\Gersten}.
His result can easily be extended to all dimensions. 
It also follows that such groups have at most
 an exponential second order \iif\ \cite{\Wang}.

We will modify Gersten's proof for groups with the \fftp, to show they are of
type \ft.
It is unclear how to extend our results to higher dimensions, 
nor to the larger class of \AC\ groups.

%%%%%%%%%%%%%%%%%%%%%%%%%%%%%%%%
\section{Fellow traveling and \ac}
%%%%%%%%%%%%%%%%%%%%%%%%%%%%%%%%
Suppose $G$ is a group with finite \gset\ $X$.
A word in $X^*$ represents a path in the \cg\ based at any vertex.
Paths can be  parameterized by non-negative
 $t\in \mathbb {R}$ by
defining $w(t)$ as the point distance $t$ along the path if $t<|w|$
and  $w(t)=\overline w$ if $t\geq|w|$, where $\overline w$ is the endpoint of $w$.
Paths $w$ and $u$ are said to {\em $k$-fellow travel} if $d(w(t), u(t)) \leq k $ for each
$ t\in \mathbb {R}$ with $t\geq 0$.
The two paths are  {\em asynchronous} $k$-fellow travelers if there is a non-decreasing
proper continuous function $\phi : [0,\infty) \rightarrow [0,\infty)$ such that 
$d(w(t), u(\phi (t))) \leq k $.
This means that any point on $w$ is within $k$ of some point on $u$ and vice versa.
We imagine the two paths traveling at different speeds (but not backtracking)
to keep within $k$ of each other.

A language $L\subset X^*$ enjoys the {\em (asynchronous) fellow traveler property} if there 
is a constant $k$
such that for each  $w,u\in L$ with $d(\overline w,\overline u) \leq  1$
in \CG, 
$w$ and $u$ (asynchronously) $k$-fellow travel.
We say a group is {\em (asynchronously) combable} if there is a language having the (asynchronous) 
fellow traveler property and surjecting to it.

\begin{figure}[ht!]
  \begin{center}
               \includegraphics[width=11cm]{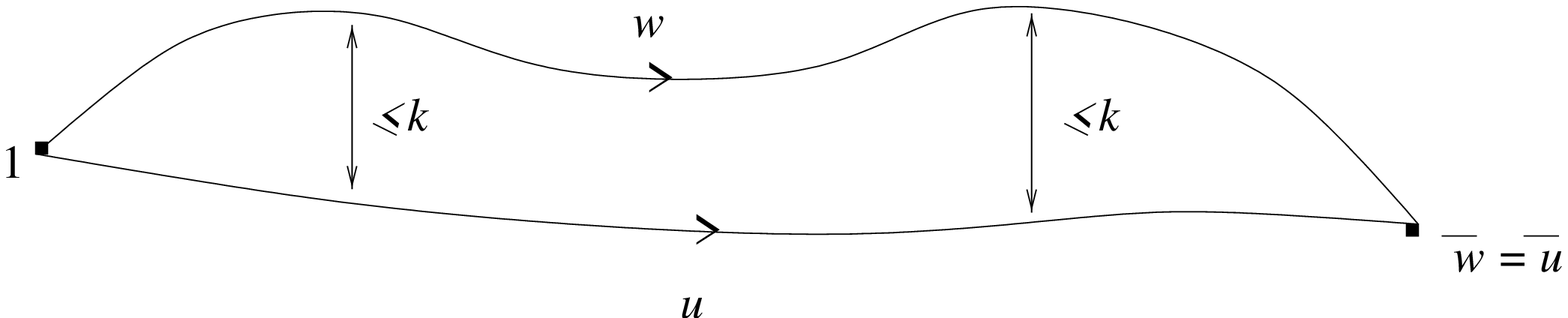}
  \end{center}
  \caption{The \fftp}
  \label{fig:1_10}
\end{figure}

A related property is the \fftp, defined as follows.
\begin{defn}
$(G,X)$ has the {\em (asynchronous) \fftp} if there exists  a constant $k$
 so that for any non-geodesic word $w\in X^*$
there exists $u\in X^*$ so that $u =_G w$, $u$ and $w$  (asynchronously) $k$-fellow travel and $|u|<|w|$ (Figure \ref{fig:1_10}).
\end{defn}

\begin{figure}[ht!]
  \begin{center}
                \includegraphics[width=11cm]{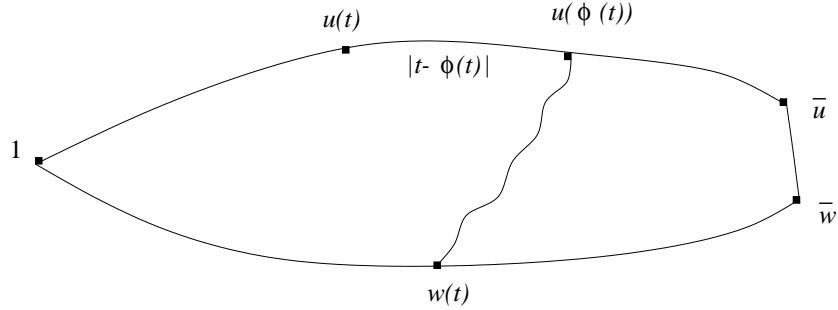}
  \end{center}
  \caption{Geodesics asynchronously $k$-fellow traveling}
  \label{fig:1_11}
\end{figure}

\begin{lem}
If $w$ and $u$ are geodesics in $(G,X)$ and asynchronously $k$-fellow travel
then they synchronously $2k$-fellow travel (Figure \ref{fig:1_11}).
\end{lem}

\noindent
Proof: 
Let $\phi:[0,\infty) \rightarrow [0,\infty)$ be a non-decreasing
proper function  such that $d(w(t), u(\phi (t))) \leq k $ for all $t\in [0,\infty)$.

\noindent
Now $u(\phi (t))$ is within $k$ of the sphere $S(t)$ of radius $t$  about 1 (namely the point $w(t)$) and $u(t)$ is a closest point of $S(t)$ to  $u(\phi (t))$, so \\
$d(u(\phi (t)), u(t)) \leq k$. Thus $d(w(t),u(t))\leq 2k$ for all $t\in [0,\infty)$.
\hfill$\Box $

\begin{cor}
The asynchronous \fftp\ and the synchronous \fftp\ are equivalent.
\end{cor}

\noindent
Proof: 
Suppose $(G,X)$ has the asynchronous \fftp\ with constant $k$.
If $w$ is not geodesic, take $t\in \mathbb N$ minimal  such that $w(t)$ is geodesic
and $w(t+1)$ is not geodesic. Let $w=w_1w_2$ with $w_1=[w(0),w(t+1)]$.
There is a word $u$ that asynchronously $k$-fellow travels $w_1$ and $|u|< t+1$.
If $u$ is not geodesic then there is a word $v$ that asynchronously $k$-fellow travels $u$
 and $|v|<|u|< t+1$. Then $v$ must be geodesic, so we have
two geodesics that asynchronously $2k$-fellow travel, so by the lemma
they synchronously $4k$-fellow travel.
Then $vw_2$ is shorter than $w$ and they synchronously $4k$-fellow travel, provided
 $k\geq 1$.
The other direction is obvious.
\hfill$\Box $

An important fact is that if $(G,X)$ has the \fftp\ then
the language of geodesics is regular \cite{\NS}.
It is not known whether  the converse is true.
%Recent work of Loeffler, Meier and Worthington \cite{Meier} provides an example that is not 
%of type \ft\ has a regular language of geodesics with 
%respect to a certain \gset. It follows from Theorem \ref{fftp} then that this implication is not reversible.

A \fg\ group is said to be {\em \AC$(i)$} if there is a constant $C(i)$ such that any two elements
of the \cg\ which lie in the metric ball  of (arbitrary) radius $n$ in the graph and lie within
 distance $i$ of each other also lie 
within distance $C(i)$ of each other in the ball of radius $n$.
Cannon proves that \AC$(i)$ implies \AC$(i+1)$
for $i \geq 2$, so we say a group is {\em \AC} if it is \AC$(2)$.

\begin{figure}[ht!]
  \begin{center}
      \includegraphics[width=10cm]{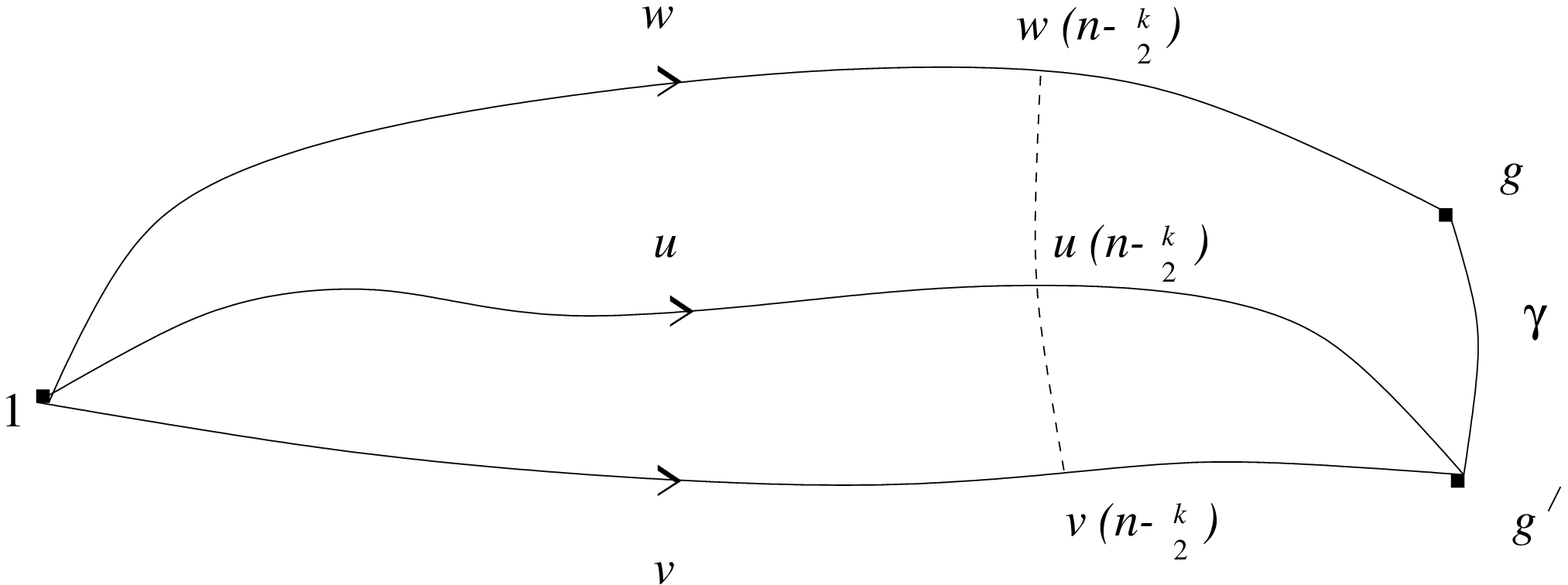}
  \end{center}
  \caption{The \fftp\ implies \ac.}
  \label{fig:4_01}
\end{figure}

\begin{propn}
If $(G,X)$ enjoys the \fftp\ then $(G,X)$ is almost convex (Figure \ref{fig:4_01}).
\end{propn}

\noindent
Proof:
Suppose $k>0$ is the \fftp\ constant,
and that $g,g' \in S(n)$ with $d(g,g')\leq 2$ realized by a path \ga, where $S(n)$ is the 
set of all points of the \cg\ that lie distance $n$ from the identity.
Let $w$ be a geodesic path for $g$. Now $w\gam$ is not geodesic for $g'$, so
by the \fftp\ there is a path $u$ for $g'$ which $k$-fellow travels $w\gam$, and
$$|u| < |w\gam| \leq |w|+ 2 =n+ 2.$$
If $u$ is not geodesic then there is a path $v$ for $g'$ which 
$k$-fellow travels $u$ and 
$$|v| < |u| < n+ 2$$
hence $v$ must be geodesic.
If $u$ is geodesic put $v= u$.

Since these paths pairwise $k$-fellow travel, then it is easily checked that
the path from $w(n-\frac{k}{2})$ to $u(n-\frac{k}{2})$ to $v(n-\frac{k}{2})$
is contained inside the ball $B(n)$ of radius $n$.
Thus we have shown that $(G,X)$ is \AC\ with constant at most $3k$.
\hfill $\Box$

An example in \cite{\Ethesis} shows that this implication is
not reversible.
We have already noted that both properties are \gset\ dependent; Thiel \cite{\Thiel}
gives an example of a group that is \AC\ for one \gset\ but not another, and Neumann and
Shapiro \cite{\NS} give a virtually abelian group which enjoys the \fftp\ 
for one \gset\ but not another.
Cannon proved that \AC\ groups are \fp\ and have at most an exponential \iif\ \cite{\Cannon}.

\begin{propn}
If $(G,X)$ enjoys the  \fftp\ then $G$ is \fp\
 and has at most a quadratic \iif.
\end{propn}

\noindent
Proof:
Suppose $w \in X^*$ is a word evaluating to 1 in $G$. Then the edge path described by $w$
in the \cg\ is a loop.
Unless $w$ is the empty word, it is not geodesic, so by the \fftp\ there is a shorter path
 which $k$-fellow travels it. Iteratively we can find successively shorter 
paths until we get the empty word. This requires at most $|w|$ iterations. Then the
 loop $w$ can be filled by at most $|w|^2$ 
relators of length at most $2k+2$, and the result follows.
\hfill$\Box $

%%%%%%%%%%%%%%%%%%%%%%%%%
\section{The main theorem}
%%%%%%%%%%%%%%%%%%%%%%%%%
One concrete way to construct an \EM\ space is to start with a presentation 
2-complex for a group, which is a space having one vertex,
 an edge for each generator and a 2-cell for each relator, glued in 
appropriately.
The fundamental group of this 2-complex is the group, and if it has nontrivial
second order homotopy we glue in 3-cells to kill it.
Inductively we can glue in higher dimensional balls to obtain an \EM\ space for
$G$.
The universal cover of this construction has the \cg\ for its 1-skeleton, and the
 {\em \cc}
(or ``filled \cg'' \cite{\Epstein}) for its 2-skeleton.

The metric on the two complex can be defined by taking the metric on the 1-skeleton as
 the metric on the \cg, and saying that a 2-cell is ``in the $n$-ball $B(n)$'' if its
 boundary is in $B(n)$ of the \cg.

%\clearpage
%\thispagestyle{empty}
\begin{figure}[p]
  \begin{center}
    \subfigure[Tops and sides matching up]{
      \includegraphics[width=8cm]{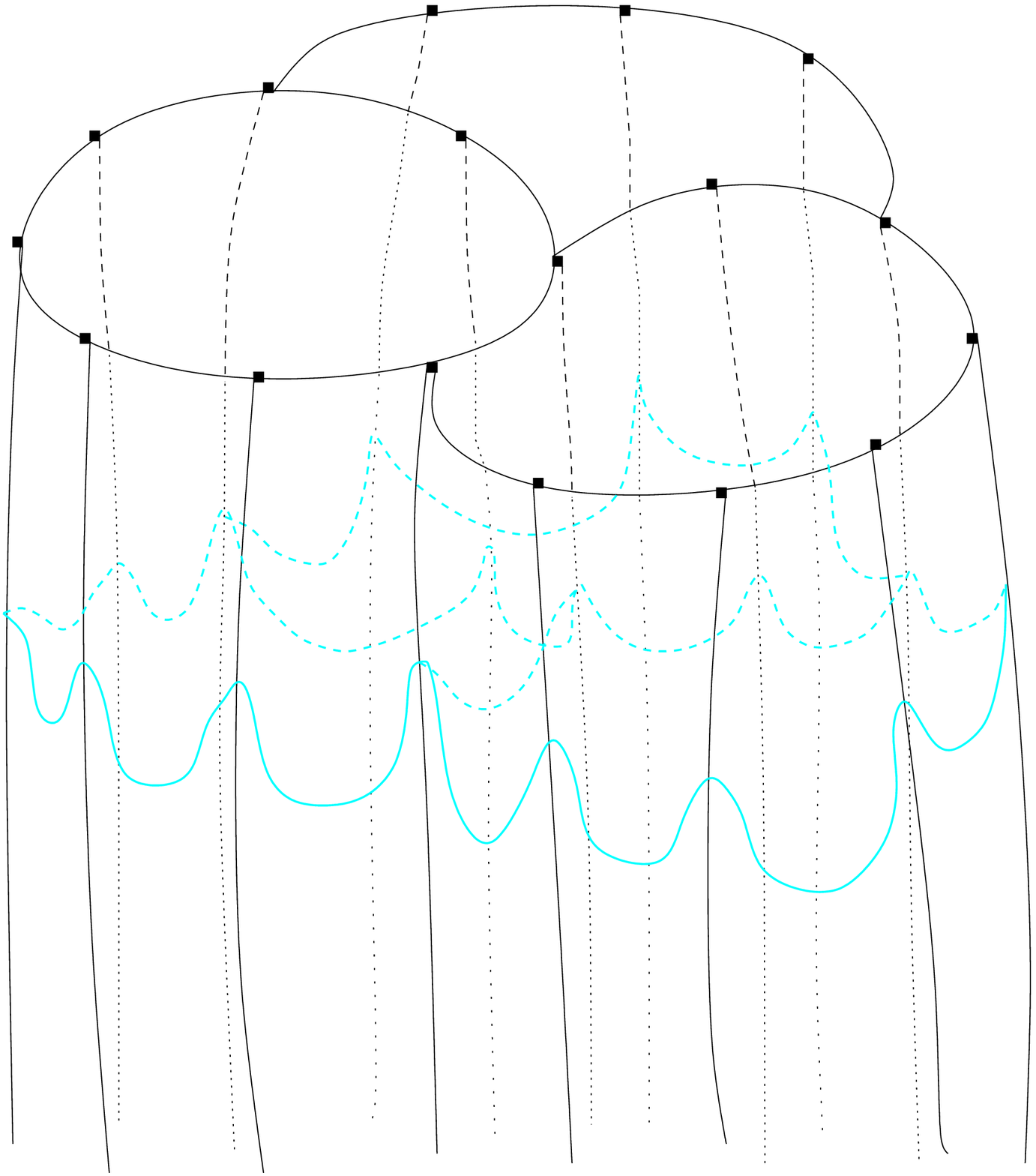}}
    %  \label{4_04_a}
      \qquad \qquad
    \subfigure[Gluing on bases]{
      \includegraphics[width=8cm]{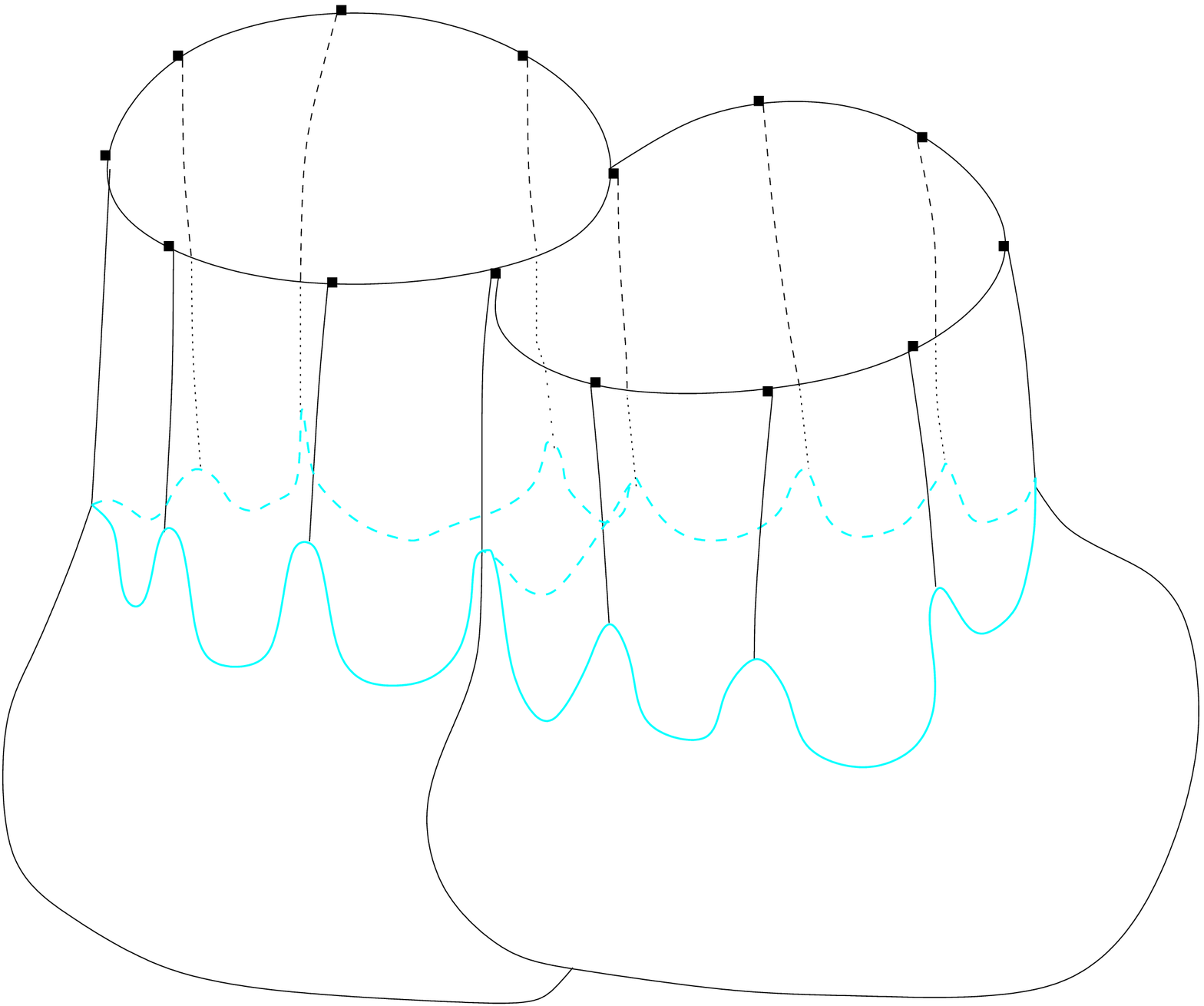}}
     % \label{4_04_b}
  \end{center}
  \caption{The drum construction}
  \label{drum}
\end{figure}

%\clearpage
The following  observation is required for the main theorem.
\begin{lem}\label{ball}
Let $G$ be any group with isoperimetric function $\rho$ and suppose
that the universal cover of an \EM\ space for $G$ has 2-cells of perimeter at most $2k+2$.
Consider a  loop $L$ of length $s$ 
with its vertices lying inside $B(r)$ (the ball of radius $r$ about the identity of $G$) in the 1-skeleton.
Then $L$ can be filled by $2$-cells so that all interior points lie in 
$B(r+(2k+2)\rho(s))$.
\end{lem}

\noindent
Proof:
The   loop $L$ can be filled by at most $\rho(s)$ 2-cells, and each cell has at
most $2k+2$ edges. Then there are are most $\rho(s)(2k+2)$
edges in this filling. A geodesic path from an interior point to a point in $L$
has no more than the number of edges in the interior, hence every point
lies in $B(r+\rho(s)(2k+2))$.
\hfill $\Box$

\begin{thm}
\label{fftp}
If $(G,X)$ has the \fftp\ then $G$ is of type \ft.
 \end{thm}

\noindent
Proof:
Let $K$ be the \cc\ for $G$ with respect to the presentation
\begin{eqnarray*}
\langle X & | & \{r\in X^*:r=_G 1, |r|\leq 2k+2\} \rangle.
\end{eqnarray*}
We will show that any combinatorial 2-sphere occurring in $K$ can be filled by 3-cells of a bounded size.
This bound will depend only on the \fftp\ constant $k$, the \iif\ $\rho $ for $G$ and 
a constant $\epsilon \in \mathbb N$.

Let $\Theta $ be a combinatorial 2-sphere of arbitrary size in $K$.
By isometry we can assume  the identity is a vertex of $\Theta$.
Let $n=\min\{m\in \mathbb N:\Theta \subseteq B(m)\}$.
For each 2-cell  of $\Theta $ we will attach a 3-cell so that the boundary of
$\Theta \cup \{$~3-balls~$\}$  is $\Theta $ and another combinatorial 2-sphere 
$\Theta '\subseteq B(n-\epsilon+\frac{k}{2})$.
Provided $\epsilon$ is chosen to be greater that $\frac{k}{2}$ we can inductively
 fill $\Theta $ by  3-cells.

Following Gersten we can think of our 3-cells as ``drums'', albeit distorted ones.
The top of each drum will be a unique 2-cell of $\Theta $.
The sides will be described presently, and they will match up with adjacent drums,
that is, drums having tops adjacent in $\Theta $.
Each drum will have a base that need not match up, as shown in Figure \ref{drum}, and we require
 that the entire base is contained in $B(n-\epsilon+\frac{k}{2})$.  
The extra $\frac{k}{2}$ here is due to some of the base cells  ``bulging out'' 
of the base (these will be ``type 2'' below).
After attaching one such 3-cell for each 2-cell of $\Theta$ the bases will glue together to form a 
homotopic copy $\Theta '$ of $\Theta$ inside
$B(n-\epsilon +\frac{k}{2})$.
This is a sketch of the argument; now let us fill in the details.

\noindent\textbf{Constructing the sides}

\begin{figure}[ht!]
  \begin{center}
      \includegraphics[width=5cm]{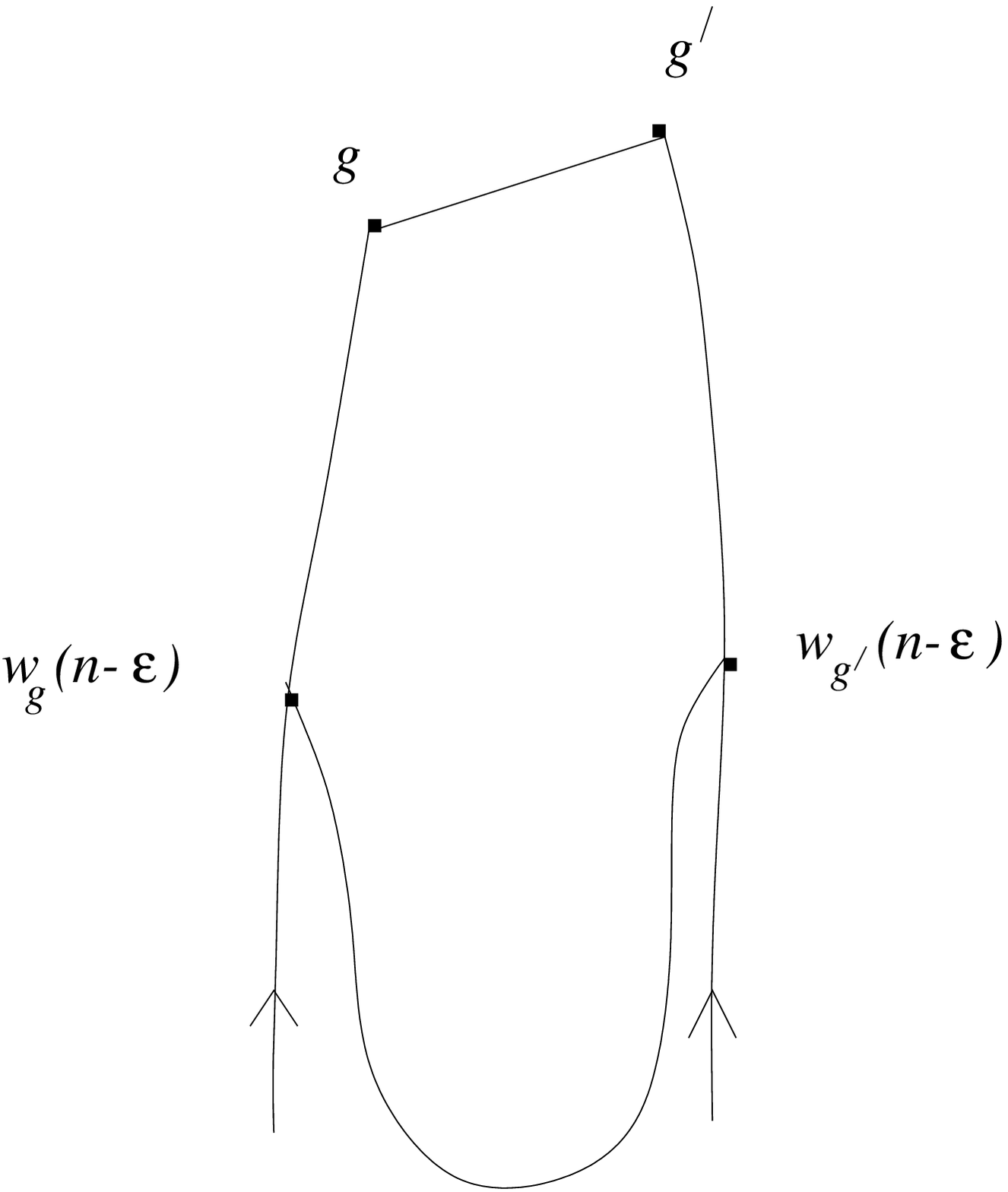}
  \end{center}
  \caption{A side cell}
  \label{fig:4_05}
\end{figure}

\noindent Fix a set of geodesics  $w_g$ from 1 to each vertex $g$ of $\Theta$.
Fix the constant $\epsilon \in \mathbb N$.
Consider each edge $(g,g')$ of $\Theta$ that lies outside of $B(n-\epsilon)$.
Retrace the geodesics $w_g,w_{g'}$ back to
 $w_g(n-\epsilon),w_{g'}(n-\epsilon)$.
Recall that since $(G,X)$ has the \fftp, $(G,X)$ is \ac\ so there is
a path from $w_g(n-\epsilon),w_{g'}(n-\epsilon)$  inside
$B(n-\epsilon)$ of length at most $C(2\epsilon +1)$
where $C(i)$ is the \ac\ constant for points distance at most $i$ apart.
Now $C(i)$ is a function of  $k$ and $i$, so
$C(2\epsilon +1)$ depends on $k$ and  $\epsilon$.
A {\em side cell} is either the edge $(g,g')$ if it lies inside $B(n-\epsilon)$, or a 
2-disc of perimeter at most $C(2\epsilon+1)+2\epsilon+1$, as seen in Figure \ref{fig:4_05}.
Each  side cell has at most $\rho(C(2\epsilon+1)+2\epsilon+1)$ 2-cells of the \cc\ $K$.
Each drum has a {\em top} a 2-cell of $\Theta $ of perimeter at most $2k+2$, so has at most
$(2k+2)\rho(C(2\epsilon+1)+2\epsilon+1)$ 2-cells for its sides.

Now we have a loop of length at most $(2k+2)C(2\epsilon+1)$ 
to which we must attach a base. This loop lies in $B(n-\epsilon)$.
Effectively we have taken the 1-skeleton of  $\Theta $ and pushed it down into 
$B(n-\epsilon)$ by homotoping each edge outside $B(n-\epsilon)$ to a path of length
at most $C(2\epsilon+1)$, the homotopy for each edge realized by a side cell.
So we have a homotopic copy of the 1-skeleton of $\Theta $ inside $B(n-\epsilon)$.
 For each 2-cell of $\Theta $ we have a loop in
this copy of length at most $(2k+2)C(2\epsilon+1)$.
We will attach a base of bounded size to each such loop, and ensure that
each base lies in $B(n-\epsilon)$.

\noindent\textbf{Constructing the base}

\noindent Each drum so far has a top, sides, and a loop of length at most $(2k+2)C(2\epsilon+1)$
in $B(n-\epsilon)$ to which we must attach a base.
Suppose the loop has length $l$. Fix a vertex on the loop, let $u_0$ be a geodesic
from 1 to it, and write the loop as an edge path
$a_1a_2\ldots a_l$. By the \fftp, if $u_0a_1$ is not geodesic then there is a
word $u_1=_Gu_0a_1$ such that $|u_1|<|u_0|+1$, so $|u_1|\leq|u_0|\leq n-\epsilon$, and
$u_0,u_1$ $k$-fellow travel.
If $u_0a_1$ is geodesic then put $u_1=u_0a_1$ and note that $|u_1|\leq n-\epsilon$ since
it is a geodesic for to point on the loop.
Recursively we can find $u_0,u_1,\ldots u_l$ such that $|u_i|\leq n-\epsilon$
and $u_{i-1},u_i$ $k$-fellow travel for $i\in [1,l]$.
Note that $u_l,u_0$ need not $k$-fellow travel, but they do $kl$-fellow travel.
We call this the ``tear'' in the drum.
\begin{figure}[ht!]
  \begin{center}
      \includegraphics[width=12cm]{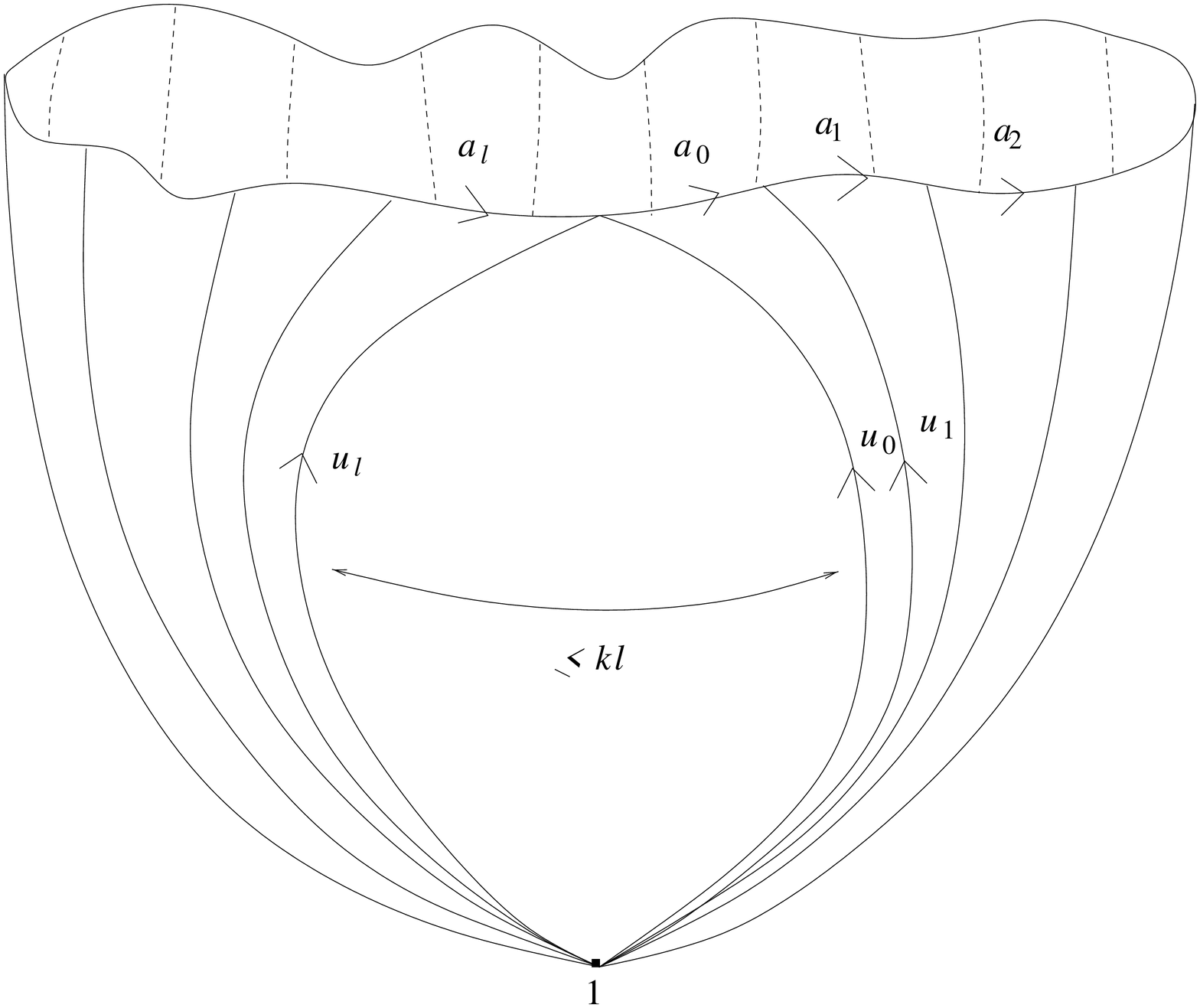}
  \end{center}
  \caption{The ``tear'' in the base}
  \label{fig:4_06}
\end{figure}

\newpage
\noindent\textbf{Type 1 base cells:}

\begin{figure}[ht!]
  \begin{center}
      \includegraphics[width=8cm]{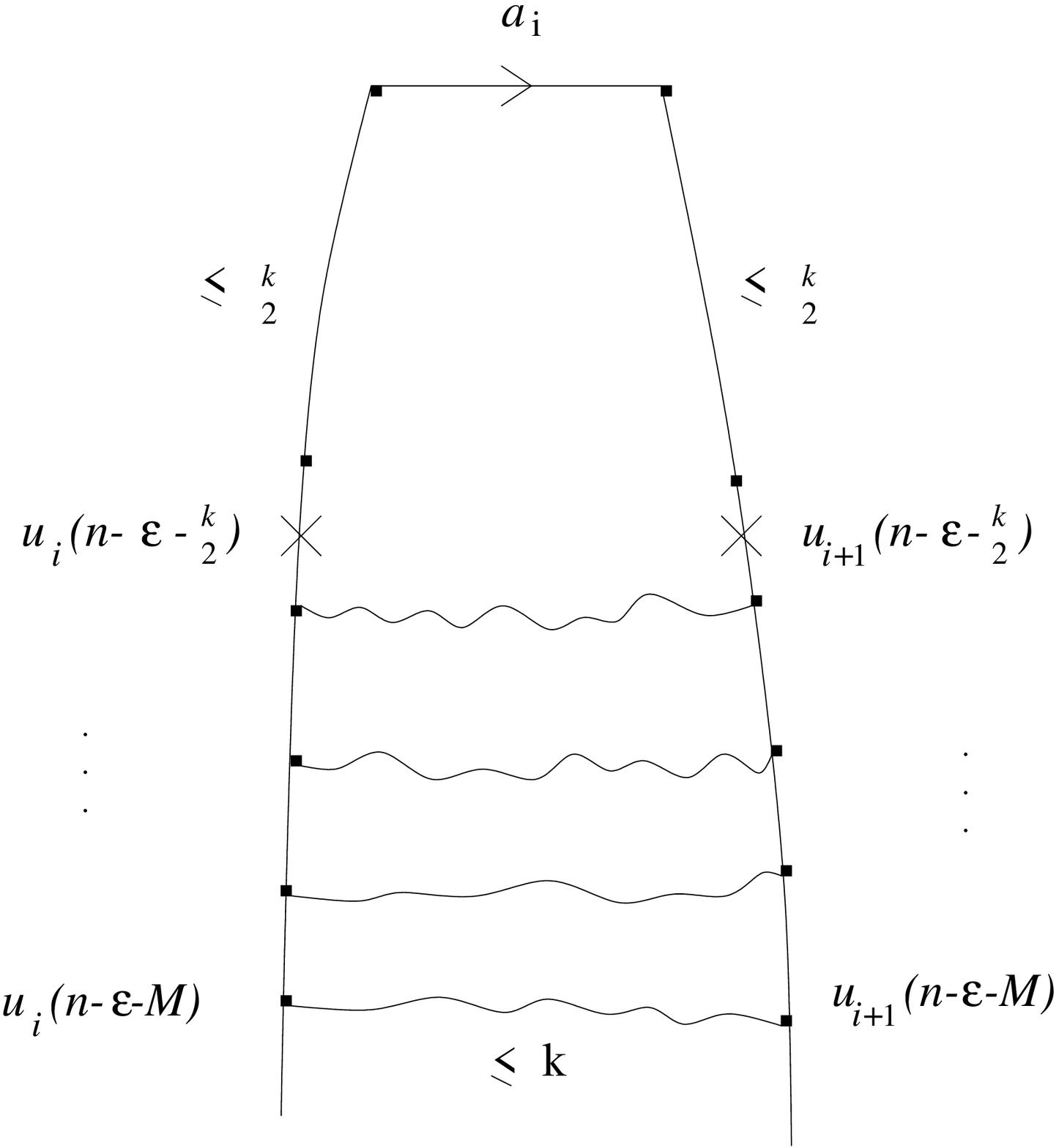}
  \end{center}
  \caption{Base type 1}
  \label{fig:4_07}
\end{figure}

\noindent Retrace each path $u_i$ back to $u_i(n-\epsilon-M)$
where $M$ is a constant to be determined below.
This gives at most $l$ type 1 cells as in Figure \ref{fig:4_07}.
The cell at the top has perimeter at most $2k+2$, and below it each cell has
 perimeter at most $2k+2$.
So in total we have at most $Ml$ 2-cells of  $K$ to make up the type 1 base cells 
for each drum.
For each integer $t\geq \frac{k}{2}$ 
 each pair of points $u_{i-1}(t), u_i(t)$ has a path of 
length at most $k$ between them, thus each of these cross paths lies in $B(n-\epsilon)$.
It follows that these cells lie in $B(n-\epsilon)$.

\noindent \textbf{Type 2 base cells:}

\begin{figure}[ht!]
  \begin{center}
      \includegraphics[width=10cm]{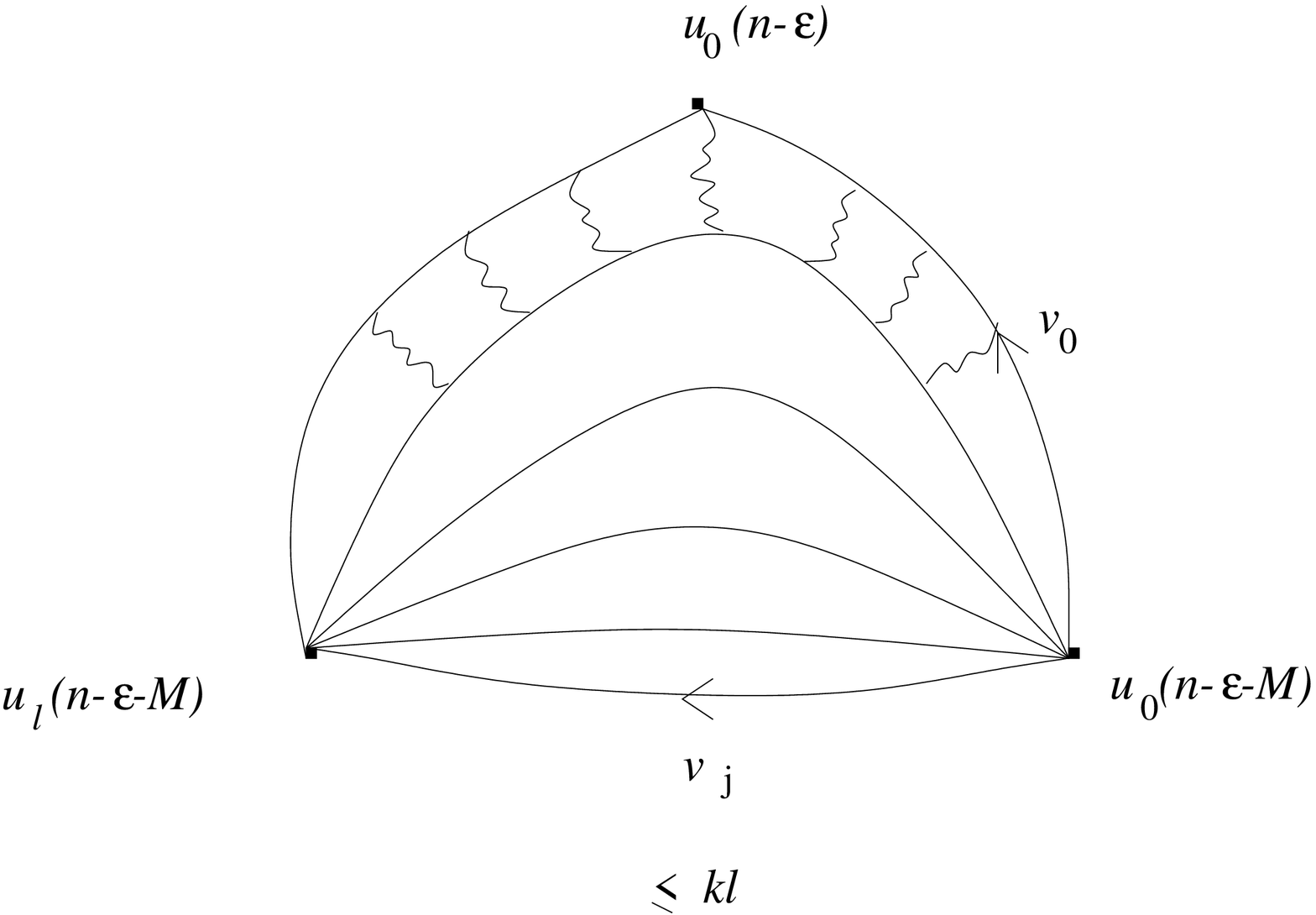}
  \end{center}
  \caption{Base type 2}
  \label{patch}
\end{figure}

\noindent 
We want to fill in the tear with cells inside $B(n-\epsilon)$. 
Let $v_0$ be the path  from 
$u_0(n-\epsilon-M)$ to $u_0(n-\epsilon)=u_l(n-\epsilon)$ to 
$u_l(n-\epsilon-M)$
(some of the points $u_0(n-\epsilon-t)$ could be the identity if $n-\epsilon < M$).
This path has length at most $2M$, and we know there is a path of length at most $kl$ for it.
So if  $v_0$ is not geodesic then there is a shorter path $v_1$ from
$u_0(n-\epsilon-M)$ to $u_l(n-\epsilon-M)$
 which $k$-fellow travels $v_0$.
Recursively if $v_i$ is not geodesic we can find a shorter path $v_{i+1}$ which $k$-fellow travels $v_i$.
This gives at most $2M-kl$ paths, as shown in Figure \ref{patch}.
Now each path $v_i$ has length at most $2M$, so lies in $B(n-\epsilon)$.
For each integer $t$ there is a path of length at most $k$ from $v_i(t)$ to $v_{i+1}(t)$, and this path lies in
$B(n-\epsilon +\frac{k}{2})$ So we can fill in the tear with at most $(2M-kl)2M$ 2-cells of perimeter at most
$2k+2$, which lie in $B(n-\epsilon +\frac{k}{2})$.
Note that we need an extra $\frac{k}{2}$  here.

\noindent \textbf{Type 3 base cells:}

\noindent After including the above base cells
in our drum we are left with a loop of length at most $2kl$ that lies in 
$B(n-\epsilon-M+\frac{kl}{2})$.
By  Lemma \ref{ball} this loop can be filled by at most $\rho(2kl)$
2-cells of $K$
so that the interior lies in
 $$B(n-\epsilon-M+\frac{kl}{2}+(2k+2)\rho(2kl)).$$
By choosing 
$$M=\frac{k(2k+2)C(2\epsilon+1)}{2}+(2k+2)\rho(2k(2k+2)C(2\epsilon+1))$$
$$\geq\frac{kl}{2}+(2k+2)\rho(2kl)$$
we ensure the type 3 base cells lie in $B(n-\epsilon)$.
Note that $\rho $ is a monotone increasing function, so these inequalities are justified.

In total each drum has a boundary of at most
\[  \left\{ \begin{array}
{r@{\quad 	 \quad}l}
 1					&  \mathrm{top}  \\
 (2k+2)\rho(C(2\epsilon+1)+2\epsilon+1) & \mathrm{sides}  \\
 (2k+2)M				 & \mathrm{base} \; \mathrm{type} \; 1 \\
 2M(2M-k(2k+2)C(2\epsilon+1))            & \mathrm{base} \; \mathrm{type} \; 2 \\
\rho(2k(2k+2)C(2\epsilon+1)) 		& \mathrm{base} \; \mathrm{type} \; 3
\end{array} \right. \]
2-cells of $K$, thus each drum has a bounded size dependent on the constants
$k,\epsilon$ and $\rho$.
\hfill $\Box $

\begin{cpfftp}
%Schollium?
If $(G,X)$ has the \fftp\ then $G$ has at most exponential second order isoperimetric
function.
\end{cpfftp}

\noindent
Proof:
Fix the 3-complex constructed above.
Suppose  a combinatorial 2-sphere $\Theta $ has area $N$, that is, it consists of  $N$ 2-cells.
By isometry we may assume $1\in \Theta $, and let $n$ be the smallest integer
such that $\Theta \subseteq B(n)$.
For each 2-cell in $\Theta $ we attach one 3-ball. Let 
$$b= (2k+2)M	+  2M(2M-k(2k+2)C(2\epsilon+1))  +   \rho(2k(2k+2)C(2\epsilon+1)) 	$$
 which is greater than the number of 2-cells
in the base of a 3-ball from the proof of the theorem above.  
Then after attaching $N$ 3-balls we obtain another combinatorial 2-sphere having area at most
$Nb$ and which lies in $B(n-\epsilon+\frac{k}{2})$.
If we repeat this procedure $\frac{n}{\epsilon-k/2}$ times we get a 2-sphere
inside $B(0)$ so it must be the identity vertex. We will have filled 
$\Theta $ with at most
$N+Nb+(Nb)b+\ldots +Nb^{\frac{n}{\epsilon-k/2}}$ 3-balls.
Now since $N$ is the area of $\Theta $, and $1\in \Theta $, then $n$ can be at most
$N(2k+2)$ which is the maximum number of edges in $\Theta $. 
Therefore the number of 3-balls required to fill a combinatorial 2-sphere of
area $N$ is
$$\sum _{i=0} ^{(2k+2).N} Nb^i=Nc^N <d^N$$ for c,d constants.
\hfill $\Box $

\subsection*{Acknowledgments}
I wish to thank my advisor Walter Neumann for his guidance and encouragement.

\bibliography{refs}
\bibliographystyle{plain}

\end{document}